\documentclass[12pt]{article}

\newtheorem{lemma}{Lemma}

\newtheorem{remark}{Remark}
\newtheorem{definition}{Definition}

\usepackage[latin1]{inputenc}
\usepackage{t1enc}
\usepackage{a4}
\usepackage{amsmath}
\usepackage{amssymb}
\usepackage{psfrag, picinpar}
\usepackage{graphicx,color}
\def\colon{\hspace{-2mm}\emph{\textbf{:}}\ }

\def\lam{\lambda}
\def\th{\theta}

\def\om{\omega}

\def\i{\infty}

\def\ra{\rightarrow}
\def\R{\mathrm{R}}
\def\Ex{\mathbb E}
\def\Pr{\mathbb P}
\def\tf{\tfrac}
\def\tes{tessellation}
\def\zwei{\R^2}
\def\lseg{line-segment}
\def\ind{independent}
\def\piv{$\pi$-vertex}
\def\pivs{$\pi$-vertices}
\def\pro{probability}
\def\calQ{\mathcal{Q}}
\def\calG{\mathcal{G}}
\def\bfG{\mathbf{G}}
\def\calX{\mathcal{X}}
\begin{document}
\baselineskip = 15pt
\title{\sc Tessellations derived from random geometric graphs}
\author{Richard Cowan\\ School of Mathematics and Statistics\\ University of Sydney, NSW, 2006, Australia\\ \textit{e-mail: richard.cowan@sydney.edu.au}\\
\\
 Albert K. L. Tsang\\ \textit{e-mail: akltsang1@gmail.com}}

\date{}
\maketitle
\begin{abstract} In this paper we consider  a random partition of the
plane into cells, the partition being based on the nodes and links
of a {\it random planar geometric graph}. The resulting structure
generalises the \emph{random \tes}\ hitherto studied in
the literature. The cells of our partition process, possibly
with holes and not necessarily closed, have a fairly
general topology summarised by a functional which is similar to the
Euler characteristic.  The functional can also be extended to certain cell-unions which can arise in applications. Vertices of all valencies, $0, 1, 2, ...$ are
allowed.  Many of the formulae from the traditional theory of random
tessellations with convex cells, are made more general to suit this
new structure.  Some motivating examples of the structure are given.
\end{abstract}
\bigskip


\medskip
\bigskip
\noindent
\section*{1. Introduction}\label{sec1}

In this paper we consider a very general type of stationary
random `\tes'\ in the plane, one which might be called a {\it random
planar partitioning} (RPP). Later in this introduction, we shall describe this general structure but, as a foundation, we firstly review the status of the random \tes\ literature.

{\sc Basic information about convex-celled \tes s:} The theory of planar random \tes s has evolved over the years partly through particular models and partly via model-free studies which assume only \emph{stationarity} and \emph{local finiteness}. Often researchers of model-free \tes s have added a convexity assumption --- cells must be convex --- and also a statement about the closed or open status of the cells; see \cite{amb70}, \cite{amb74}, \cite{cow78}, \cite{cow80}, \cite{mec80}, \cite{mec84} and \cite{skm95}.

We commence our discussion with the version where all cells in the model-free \tes
\ are \emph{closed, bounded} and \emph{convex}. Firstly we define a \tes\ in this context.

\begin{definition}\hspace{-2mm}\emph{\textbf{:}}\label{defn_tess}
A \textbf{\tes}\  of the plane is a locally-finite collection of compact convex
 \textbf{cells}, each of positive area, which cover the plane and overlap only on
cell boundaries\footnote{ Convexity of cells implies, of course, that cells are polygons. A collection is locally finite if
every bounded domain in $\zwei$ intersects a finite number of
cells.}. The union of the cell boundaries is called the
\tes\ \textbf{frame}. Each cell, being a polygon, has
\textbf{sides}\footnote{The cell's sides and other \lseg s are
considered closed sets, unless stated otherwise.}
 and \textbf{corners}; they lie on the frame. The union (taken
 over all cells) of
 cell corners is a collection
 of points in the plane called
 the \textbf{vertices} of the \tes. Those closed line-segments which are
 contained in the frame, have a vertex
 at each end and no vertices in their relative-interior are called \textbf{edges} of
 the \tes.
The number of edges emanating from a vertex is called the
\textbf{valency} of that vertex.
\end{definition}

 Whereas no vertices can lie in the interior of an edge, some vertices might
 lie in the interior of a cell-side.

\begin{definition}\hspace{-2mm}\emph{\textbf{:}}\label{defn_pivs}
    Vertices which lie in the relative interior\footnote{The terminology `relative interior' is technically more
    correct for \lseg s imbedded in the plane. We shall usually, however,
    drop the word `relative' in the rest of the paper. `Interior' means `relative interior'.} of a cell-side
    are called \textbf{\pivs}. The name arises because in such vertices one of the angles between
    consecutive emanating edges is $\pi$. 
\end{definition}

The
focus of attention in many studies of planar \tes s and tilings has
been the
 \emph{side-to-side} case.

\begin{definition}\hspace{-2mm}\emph{\textbf{:}}\label{defn_sts}
    A \tes\ is \textbf{side-to-side} if each side of any polygonal cell in
    the \tes\
coincides with a side of another cell. Alternatively, we say that
a \tes\ is side-to-side if it has no \pivs.
\end{definition}

{\sc Introducing randomness, stationarity and ergodicity:} Let $\Omega$ be the space of \tes s (each $\om \in \Omega$ conforming to
Definition \ref{defn_tess}). A convex-celled \emph{random} \tes\ is a randomly selected (or randomly
constructed) entity $\om \in \Omega$. More precisely, we place a
\pro\ measure $\Pr$ on a suitably-large class $\calQ$ of subsets of $\Omega$. So the triple $(\Omega, \calQ ,\Pr )$ is our random convex-celled \tes.\footnote{Schneider and Weil \cite[p.19]{scwe08} have called this approach, which identifies a geometric entity and a random element, the \emph{canonical} representation. Although disadvantages of this method emerge as theories become more elaborate, we use it in this paper because (a) it allows the reader to visualise a random element $\om$, (b) it more comfortably and directly connects to the concepts of ergodic theory and (c) our theory in this paper doesn't encounter those disadvantages.} The
expectation associated with $\Pr$ is written as $\Ex(\cdot)$.

  For $x,\ t \in \zwei$ define the
translation operator $T^t: x \rightarrow x + t$.  Thus $T^t$, which is also defined on
$\Omega$, translates any realisation $\omega$ by $t$.  We assume that the random \tes\
process is stationary and ergodic, via the assumption that $T^t$ (when
defined on $ \calQ$) is a measure-preserving and ergodic operator.
Intuitively,
`measure-preserving' means that  the statistical properties of the  structure are
invariant under translation. Ergodicity has implications
when taking large-domain spatial averages in the random process; it implies that spatial averages (for example, the average vertex valency) and various proportions (for example, the proportion of vertices which are \pivs) calculated inside the ball $B_r$, of radius $r$ and with centre at the planar origin, converge with probability one to a constant as $r \ra \i$.

There are three basic types of ergodic \tes. In one type, the \emph{periodic} \tes s, there is a
sub-collection of cells which forms a repeating structure; the full \tes\
covering the whole plane is made up of suitably translated copies of
this structure. It is clear that spatial averaging applies; for
example, the large-domain limit of average cell area will obviously
converge to the average cell area within the `repeating sub-collection'.
To make such a periodic \tes\ stationary, one places it on the plane
so that the planar origin is uniformly distributed within one copy of the repeating
sub-collection.

In the second ergodic type, the \emph{mixing} \tes s, the method of constructing the stationary random
\tes\ is such that features which are a considerable distance apart
are effectively statistically independent. The geometry of a \tes\
will, of course, impose a short-range dependence but this decreases
with distance in this type of ergodic \tes, in the limit (as
distance tends to infinity) to complete independence. As is well known, averaging over
independent entities leads to almost-sure convergence to a constant. For mixing \tes s
the short-range dependencies are dominated by the vastly greater
number of long-range `independencies' --- and convergence to a constant still occurs.

As for the third type, these are combinations of the first two --- for example, a periodic \tes\ modified by random operations which have a tendency toward independence as distance increases. There are many ways to combine the notions of periodicity and mixing whilst retaining the \emph{large-domain limiting condition} that comes with ergodicity.

{\sc Relaxing the convexity assumption: } One can generalise slightly the \tes\ that we have defined in Definition \ref{defn_tess}, relaxing convexity of the cells and
permitting vertices of valency $2$, but retaining \lseg\ edges. So cells remain as simple polygons, though not necessarily convex.  Such model-free \tes s have been studied in \cite{mil88}, \cite{sto86} and \cite{cots94}.

The tiling literature also allows non-convex cells (see Gr\"{u}nbaum and Shephard \cite{grsh87}). The cells, assumed closed, may have fairly general shapes that are isomorphic to a closed planar disk. Whilst this allows considerable freedom in the type of non-convex cell, the other regularity conditions used by the tilers (N.2 and N.3 imposed on page 121 in \cite{grsh87}) are far too restrictive for us, especially with our emphasis on random \tes s. Gr\"{u}nbaum and Shephard were mainly concerned with non-random tilings; N.2 and N.3 are highly appropriate regularity conditions for non-random tilings.

Another style of generalisation  is due to Z\"ahle and co-authors Weiss and Leistritz (\cite{weza88}, \cite{zah88} and \cite{leza92}). In the most recent of these studies, the edges and cells are simply-connected compact submanifolds (of dimensions $1$ and $2$ respectively imbedded in $\zwei$) with boundary. Various rules govern how cells, edges and vertices interconnect. These rules, listed below, come from the theory of $d$-dimensional \emph{cell-complexes} with $d=2$.
\begin{enumerate}
  \item The intersection of two cells is contained in the boundary of each of these cells; it is \emph{either} empty \emph{or} it is an edge \emph{or} it is an isolated vertex.
  \item The intersection of two edges is \emph{either} empty \emph{or} it is a vertex which lies at a terminus of each of these edges.
  \item Any edge is contained in the boundary of some cell.
   \item Any vertex is located at a terminus of some edge.
  \item The boundary of any cell is the finite union of some edges.
  \item The boundary of an edge is the union of two vertices.
\end{enumerate}

So the edges can be curved in the Z\"ahle/Weiss/Leistritz theory, possibly with discontinuities of slope, provided no vertex is positioned at these `corners' of the curve. Thus  vertices of valency $2$ are not allowed and, as we shall see, much of the complexity that we introduce in the next section is not allowed in their theory.

\section*{2. Generalising the planar graph }

In order to understand other natural models for partitioning the plane, we have recognised the need to generalise further than can be achieved by the techniques described above. We do this by casting the discussion in terms of planar graphs.

The frame of a \tes\ complying with Definition \ref{defn_tess} can be viewed as an infinite planar graph. This graph has some imposed geometry and some topological constraints.  As is well known, a planar graph has \emph{nodes} placed in the plane, with \emph{links} connecting some pairs of distinct nodes. In a graph which is the frame of a `Definition \ref{defn_tess} \tes', the \tes\ vertices (which must have valency $\geq 3$) play the role of the graph's nodes whilst the \tes\ edges are the links, assumed non-directed. The links must be \lseg s whose relative interiors are disjoint. The polygonal cells of the \tes\ play the role of the graph's \emph{faces}; so these faces are convex, closed and bounded. Terminology such as `node valency' and `$\pi$-node' become defined on the graph via their \tes\ meaning.

We now allow the infinite planar graph to be much more general, whilst retaining the rule that \emph{each link is a \lseg\ which does not intersect other links except at its terminating nodes}. Firstly a countable collection of distinct points in $\zwei$ are identified as the nodes. Then a countable collection of links are added;
each link, which is assumed to be a \lseg\ with a node at each end, does not intersect any other link except at these nodes.
\begin{remark}\colon
    Because edges are line-segments, there can be no more than one link between any pair of nodes. Also \textbf{loops}, a curved link from a node to itself, cannot occur.
\end{remark}

\begin{figure}[h]
    \psfrag{c}{$\pi$}
    \begin{center}
        \includegraphics[width=120mm]{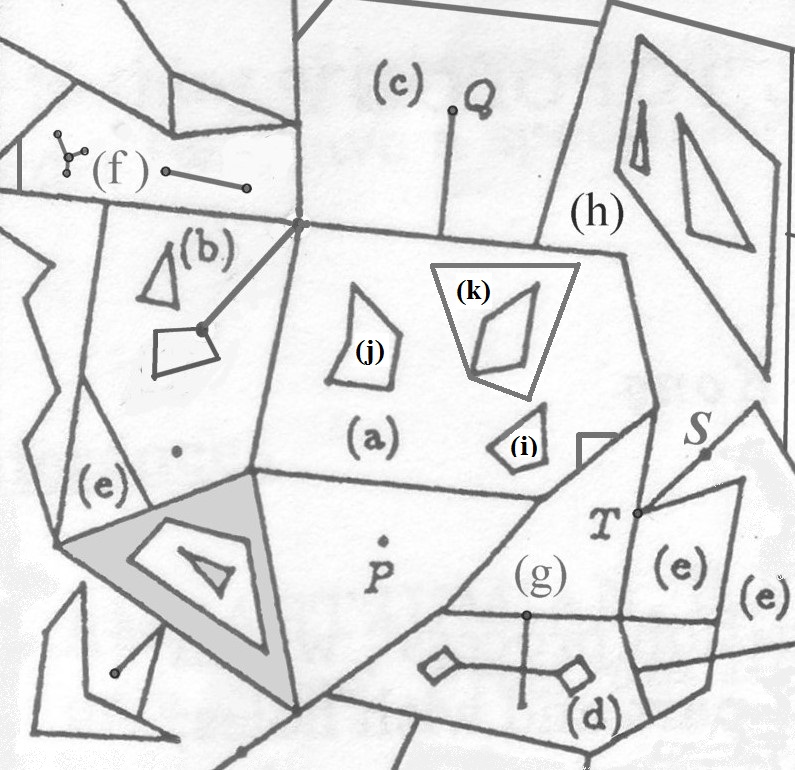}
 \\
        \caption{\label{fig1} \scriptsize  An infinite geometric graph, as seen within a window. It can also be viewed as a new sort of \tes\ having fairly general cells and some status given to a few cell-unions. Using the \emph{node, link} and \emph{face} terminology of graph theory, we see that
$P$ is a node of valency $0$, whilst $\calQ $ is a node of valency $1$. There are many $2$-valent nodes, mostly with non-collinear links emanating; node $S$ is a $2$-valent example with collinear links; it is called a double-$\pi$ node. Node $T$ is one of the many $\pi$-nodes. Viewing the figure as a \tes, the three cells marked (e) form a cell-union and the two cells shaded a darker grey are the parts of a cell-union; such unions would need to be based on a natural nexus (though the reader won't see any natural nexus in these two cases). The three cells (i), (j) and (k) have a natural nexus, however, each being cellular holes in the same cell (namely (a)). So we might like to consider (i), (j) and (k) as a cell-union.  The text below explains why cell-unions are studied.
Cell (a), with three holes, has Euler Entity $-2$. Cell (f) has two holes (each of area zero) comprising three edges and one edge respectively. Other features are discussed in the text.}
    \end{center}
\end{figure}

Figure \ref{fig1} illustrates the graph structure; in particular, the figure shows how complicated the graph's faces have become. Indeed the definition of a \emph{face} now requires considerable care.

{\sc The concept of a face:} In the following definition, the boundary of a set $F\subset \zwei$ is defined to be $cl(F)\,\cap\,cl(\bar F)$, where $cl(\cdot)$ indicates closure and $\bar F$ is the complement of $F$. The interior of $F$ is written $int(F)$.
\begin{definition}\hspace{-2mm}\emph{\textbf{:}}\label{defn_face}
    Let $\calG$ be our infinite geometric planar graph defined by its collection of nodes and links. Let $\calG_\cup$, called the graph-union, be the subset of $\zwei$ defined as the union of all nodes and links. An open subset $F$ of $\zwei$ is called a \textbf{face} of the graph $\calG$ if and only if $F$ is connected, $F \cap \calG_\cup = \emptyset$ and the boundary of $F$ is contained in $\calG_\cup$. Subsets $F$ which are not open are not faces.
\end{definition}
\textbf{Assumption:} \emph{We assume that $\calG$ is such that $\calG_\cup$ has no \emph{unbounded} faces. We also assume that $\calG_\cup$ is \emph{locally finite}: that is, every bounded domain of $\zwei$ intersects a finite number of faces.}

Our decision to treat each face as an open set deserves comment, because the reader will notice that the occasional face --- for example, (a), (g), (h), (i) and (j) in Figure \ref{fig1} --- can be considered as a closed set\footnote{Or as a set which is neither open nor closed!} without changing the graph-union in any way. Many faces, however, have a closure which destroys part of the face's boundary (and part of the graph-union too) and this significantly alters the face's topological status. See for example faces (b), (c), (d), (f), (k) and the open `quadrilateral with a hole' that surrounds the isolated node $P$. Also, the open faces (d) and (k) are  simply-connected, but their closures are not -- and open cell (b) has one hole, its closure two. So certain features of faces naturally occurring in a graph-union would disappear if cells were considered as closed sets. It is important therefore to have a theory based on \emph{open} faces; thereby greater diversity in the graph-union occurs.


\begin{remark}\hspace{-2mm}\emph{\textbf{:}}
    Modern definitions of a graph face when the graph $\calG$ is planar and infinite, but not connected, are rare. Beineke's \cite{bei97} definition, which wins for brevity, is as follows. \texttt{Maximal connected sets in the planar set $\zwei \setminus \calG_\cup$ are called faces.} We have not adopted his definition. It would be equivalent to Definition \ref{defn_face}, however,  if it commenced \texttt{Maximal connected open sets ...}.
\end{remark}

{\sc The generalised \tes:} We started this paper with a traditional convex-celled \tes\ of the plane, then transferred our thoughts to its associated geometric planar graph. Then we allowed the graph to have more of the features that geometric graphs can have. Now, as suggested above, we reverse the transfer and look at the planar \tes--like structure generated by the more elaborate graph. Words such as `node', `link', `face' and `graph-union'  return to the more familiar `vertex', `edge', `cell' and `frame' in a \tes\ context. Whilst we might still call this structure a \tes\ (albeit described as a \emph{generalised \tes}\ or a \emph{\tes\ derived from a geometric graph}), we also call it a \emph{planar partitioning}.

Both the graph $\calG$ and its derived \emph{generalised \tes}\ have a new concept not evident in the convex-celled theory: the \emph{double-$\pi$ node or double-$\pi$ vertex}.
\begin{definition}\hspace{-2mm}\emph{\textbf{:}}\label{double}
    A node in $\calG$ is called a \textbf{double-$\pi$ node} if, when marked on $\calG_\cup$, it is $2$-valent with collinear edges emanating.  In other words, a node is a double-$\pi$ node if and only if it is a $\pi$-node of valency two. Here, a \textbf{$\pi$-node} is a node with at least one angle formed by consecutive emanating links equal to $\pi$ --- and if the valency is two, there are two such angles, so the node is `double-$\pi$. A vertex in the derived \tes\ is a \textbf{double-$\pi$ vertex} if, in the graph context, it is a double-$\pi$ node --- and it is a \textbf{\piv}\ if it is a $\pi$-node.
   \end{definition}
Node $S$ in Figure 1 is an example of a double-$\pi$ node; it is, of course, a double-$\pi$ vertex in the \tes\ context.

\begin{remark}\colon
Does the graph-union $\calG_\cup$ contain all the information that the graph $\calG$ has? No, not unless we make sure that the double-$\pi$ vertices are specially marked, as mentioned in our phrasing of Definition \ref{double}. Without marking, these vertices are visually lost in drawings. Therefore, we make a special notation $\calG_\cup^+$ to indicate the `marked infinite graph-union': namely $\calG_\cup$ with all the double-$\pi$ vertices marked. Although the information in $\calG$ is the same as that in $\calG_\cup^+$, we refer to $\calG$ as the \emph{graph} and $\calG_\cup^+$ as the \emph{generalised \tes} (or \emph{planar partitioning} or \emph{\tes\ derived from $\calG$}). Any statement in the sequel for $\calG$ holds also for $\calG_\cup^+$, but not necessarily for the unmarked $\calG_\cup$.
\end{remark}

\begin{definition}\hspace{-2mm}\emph{\textbf{:}}\label{multipart}
    A \textbf{cell-union} of a \emph{\tes\ derived from a graph}  is a \emph{finite union} of some cells of the \tes.
\end{definition}
For example, we have marked a three-celled cell-union (e) in Figure \ref{fig1}. We note again the usefulness of an open-cell theory; the three cells involved in (e) are assumed open sets. Their union comprises three cells, each being a connected set; if we treated cells as closed sets, the union would comprises only two connected parts.

{\sc Discussion:} By this process of generalisation, we allow planar partitionings which
have disconnected features (see Figure \ref{fig1} and its caption).  A cell, though still bounded and connected, might not be simply-connected (that is, it might
have `holes').  The frame might not be connected;  this will be the case if a cell
has another cell or cluster of cells wholly enclosed within its interior; the edges
of the enclosed cell(s) will be disconnected from most other edges of the
graph.  Additionally we allow the existence of vertices of valency 1 or 0,
the latter type being simply isolated points.  The edges are closed line-segments, however,
as before.

Perhaps most importantly there are many violations of the rules used by Z\"ahle \emph{et. al}. All cells are open, therefore not compact.  The vertex $P$ contradicts their Rule 4.  Cell (a) is not simply-connected, and so on! In short, the rules of Z\"ahle \emph{et al}, when still meaningful with cells so general, are often violated. Put simply, our planar-partitioning structures are not \emph{cell complexes}.


\section*{3. Counting cell sides, corners, edges and vertices}

{\sc Sides and corners of cells:} There is a need to define a \emph{side} and a \emph{corner} of these unusual cells. Our definition involves the concept of a \emph{walk} on the graph $\calG$.
\begin{definition}\hspace{-2mm}\emph{\textbf{:}}\label{defn_walk}
     Consider a sequence of $n\geq3$ nodes from $\calG$ such that consecutive nodes in the sequence have a link between them. A \textbf{walk} on $\calG$ is such a sequence beginning and ending with the same node (which we call the walk's {\sc home}), without containing {\sc home} again in the sequence.
\end{definition}
For example, if we have nodes labelled $1,2,3$ and $4$ with non-directional links $\{1,3\}$, $\{2,3\},\{2,4\}$ and $\{3,4\}$, then the sequence $w_1 = (1,3,2,4,3,2,3,1)$ is a walk whose {\sc home} is $1$.

 So a walk contains its {\sc home} node exactly twice and may contain the other nodes in the sequence more than once. A walk may equivalently be thought of as a journey on $\calG_\cup^+$, visiting the nodes (and the implied \emph{connecting links}) in the order given by the sequence --- a journey that always returns to its starting node, {\sc home}. In the example, note that $w_2 = (4,3,2,3,1,3,2,4)$ is a walk, different from $w_1$, despite $w_1$ and $w_2$ having  journeys that visit the same nodes in the same `cyclic order'.

\begin{definition}\colon
    A \textbf{first-exit walk} is defined as a walk which `exits' each node visited (except {\sc home}) on the link which gives the walker the maximum anti-clockwise turn of his body --- but if no link involves an anti-clockwise turn, he makes the minimum clockwise turn. If the node is of valency $1$, then the walker makes a clockwise turn of $\pi$ and exits the node back along his entry link. \textbf{The turning angle} is denoted by $\zeta$ and it lies in the range $[-\pi, \pi)$ where anti-clockwise is deemed positive and clockwise negative. An angle $\zeta = 0$ applies if the walker doesn't turn at all. At the conclusion of a first-exit walk, returning to {\sc home}, it is assumed that the walker turns to face his starting direction. So this last turning angle is assumed to be part of a first-exit walk.
\end{definition}

For example entering node $T$ from above, the walker exits along the edge leading to node $S$. Since $S$ has valency $2$, its exit is by the `straight-ahead' link (the only link available). Approaching the $1$-valent node $Q$ from below, the first exit is back along the link of entry, so $\zeta=-\pi$.

\begin{definition}\colon\label{defn_facecircuit}
    Let $F$ be a face of $\calG$; by assumption $F$ is bounded. A first-exit walk where every node and link in the walk's sequence lies on the boundary of $F$ and where, when traversing every link of the walk, there is always an open neighbourhood of the walker left of the link and contained in the interior of $F$, is called a \textbf{face-circuit}. There may be more than one face-circuit of the face $F$. The \textbf{link-count of a face-circuit} is the number of link-traversals (so a link traversed twice scores $2$). The \textbf{node-count of a face-circuit} is the number of node-visits made in the face-circuit, counting node {\sc home} only once. A face-circuit also has a \textbf{corner-count} defined as the number of direction changes in the face-circuit (that is, the number of non-zero turning angles $\zeta$, the nodes in the face-circuit where $\zeta \ne 0$ being called \textbf{corners of the face circuit}). The \lseg s in the face-circuit between consecutive corners are called \textbf{sides of the face circuit}; so the face-circuit also has a \textbf{side-count}.
\end{definition}

These definitions, defined above for face-circuits, apply also to \emph{faces}.

\begin{definition}\colon\label{defn_facecounts}
   The \textbf{corners of a face} $F$ are the corners on all $F$'s face-circuits, so the \textbf{corner-count of a face} is the sum of the corner-counts of all face-circuits. Likewise for \textbf{sides of a face} and \textbf{side-counts of a face} and also \textbf{link-counts of a face}. The \textbf{node-count of a face}, however, is the sum of the node-counts for the component face-circuits \texttt{plus} the number of $0$-valent nodes that form holes in the face.
\end{definition}

For example, face (a) has four face-circuits with link-counts $10, 5, 4$ and $4$ and side-counts $8, 5, 4$ and $4$, so face (a) itself has link-count $23$ and side-count $21$. Face (b) has two face-circuits with link-counts $11$ and $3$, so face (b) has link-count $14$. Faces (c) and (d) each have just one face-circuit with link-counts $9$ and $24$ respectively. Face (f) has an `outer' face-circuit with $7$ link-counts (but only $6$ side-counts) and two `inner' face-circuits with $2$ and $6$ link-counts.

Note that for three of the four face-circuits of (a), the travel direction of the circuit is clockwise (as face (a) must be on the left).  For any of the three faces (i), (j) and (k) which make a hole in (a), the face-circuits are travelled anti-clockwise (keeping the `hole-cell' to the left).

Clearly a face-circuit's node-count always equals its link-count. Its side-count always equals its corner-count.

The word `cell' replaces 'face' when our discussion turns to generalised \tes s.
\begin{definition}\hspace{-2mm}\emph{\textbf{:}}\label{defn_edgecount}
     In the generalised \tes\ induced by $\calG$, a \textbf{cell} is equivalent to a face of $\calG$ and a \textbf{cell-circuit} is equivalent to a face-circuit. So the entities  \textbf{edge-count of a cell}, \textbf{side-count of a cell}, \textbf{corner-count of a cell} are essentially defined in Definitions \ref{defn_facecircuit} and \ref{defn_facecounts}. \textbf{The vertex-count of a cell} follows the definition of the node-count of a face in those definitions.
\end{definition}
\begin{remark}\colon
    A concept of a $\pi$-vertex (and double-$\pi$ vertex) can be defined using \emph{face-sides}. A vertex that lies in at least one face-side interior is called a \piv. A vertex that lies in two face-side interiors is called a double-$\pi$ vertex.
\end{remark}

\section*{4. Descriptor of the cell's topology}\label{sec4}

The topology of a cell is summarised by a functional rather like the \emph{Euler Characteristic} $\chi$,
defined loosely as the number of parts minus the number of `holes'.\footnote{For this purpose, a
0-vertex creates a hole in the cell which surrounds it, as do isolated edges and their end vertices, as seen in cell (f).} We call this functional, which we define in this section, by different terminology: the \emph{Euler Entity}. We use a different name because some readers of our theory have remarked that the Euler Characteristic is not usually defined on open sets\footnote{When drafting this paper we shared the view of these readers, because we were unaware of Groemer's early work \cite{gro72}, where he extended the Euler Characteristic to finite unions of polygon-interiors. Indeed his work is in $d$ dimensions, extending the Euler Characteristic to finite unions of polytope-interiors. If we adopt Groemer's definitions, much of the discussion in the next sub-section becomes redundant and our \emph{Euler Entity} is equivalent to the \emph{Euler Characteristic}.} --- and our theory produces cells which are open sets.   So the
statistical properties of the `typical' cell include mean values of
topological features such as $\chi$ and also, of course, geometric features such as area
and perimeter plus various combinatorial entities.

{\sc  Introduction of $\chi$:} There are many different contexts in the topological literature where the Euler Characteristic is defined for a set $F \subset \zwei$. Mostly, for a valid definition, the set needs to be \emph{closed}; for example, in some theories $F$ should be in the convex ring.\footnote{The convex ring comprises all finite unions of compact convex sets.} Even in the Gauss-Bonnet context, the usual discourse assumes that $F$ contains its boundary. Our method is essentially of Gauss-Bonnet style, but $F$ is now \emph{open}; so $F$ doesn't contain (or even intersect) its boundary.

 \begin{figure}[h]
    \psfrag{c}{$\pi$}
    \begin{center}
        \includegraphics[width=120mm]{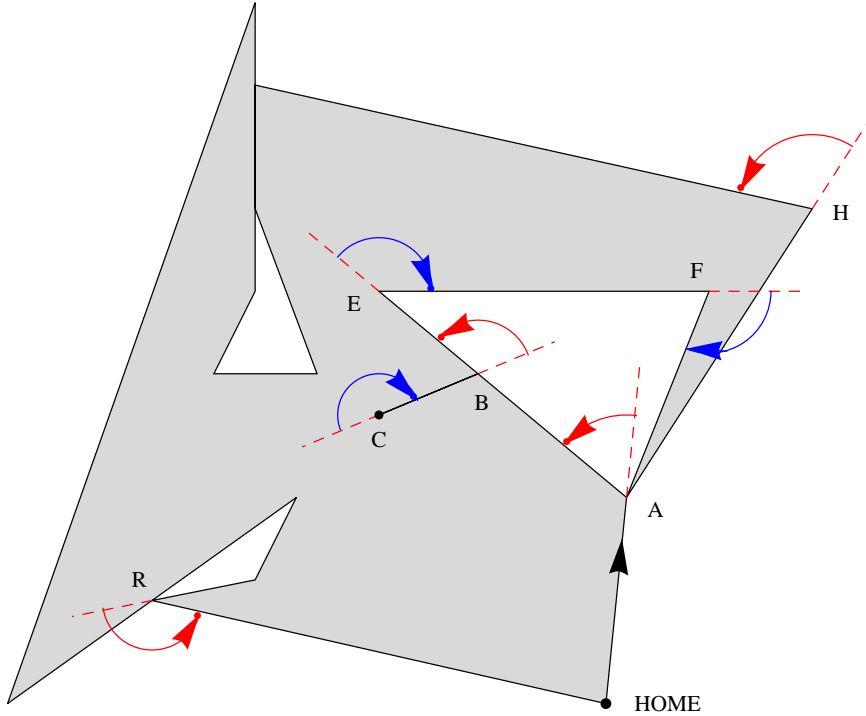}\\
        \caption{\label{fig2} \scriptsize  The grey domain here is an open face with a side-count of $19$. Its boundary is shown in black. We begin a walk on this boundary at {\sc home}, visiting nodes $A, B, C, B$ again, $E, F, A$ again, $H, ...$,  keeping the `grey face at the left' rule as links are traversed. Turning angles at some of the nodes visited --- $A, C, B, E, F, H$ and $R$ --- are shown by red or blue arrows, blue indicating a negative turning angle.  Corner $A$ is visited twice but we only show the first visit's turning angle. Only the arrow of the second visit to $B$ and to $R$ is shown.}
    \end{center}
\end{figure}

We base the idea on the `turning angles' in $F$'s face-circuits (see Figure \ref{fig2}). As explained above, at each node visit of a face-circuit the walker turns through an anti-clockwise angle $\zeta \in [-\pi, \pi)$ before exiting the node along the `first-exit link'. An angle $\zeta=0$ applies if the walker doesn't turn at all and, in general, $\zeta$ is measured from this `collinear entry and exit edges' situation. By assumption, an anti-clockwise turn gives a positive angle $\zeta \in (0,\pi)$ whilst a clockwise turn (which only occurs if no exit link involving an anti-clockwise turn is available) yields a negative $\zeta \in [-\pi, 0)$.

For example when $F$ in Figure \ref{fig1} is face (h), the walker arriving at $T$ from above turns an angle $\zeta$  of approximately $+140^\circ$ (to now walk toward $S$). When $F$ is face (c), arriving at $Q$ from below the walker turns clockwise $180^\circ$, so the turning angle is $\zeta = -\pi$. For an arrival at $S$, $\zeta =0$. When $F$ is face (g), an arrival at $T$ from below has $\zeta =0$.

The following lemma, illustrated in Figure \ref{fig2}, is trivially true.
\begin{lemma}\hspace{-2mm}\emph{\textbf{:}}\label{turn}
    For a face-circuit in our geometric graph $\calG$ (or for a cell-circuit in the planar partitioning  $\calG_\cup^+$), the sum of all turning angles is $2\pi$ if the circuit is anti-clockwise, as in Figure \ref{fig2}. This sum is $-2\pi$ if the circuit is clockwise.\hfill$\square$
\end{lemma}

Thus we are led, in the spirit of the Gauss-Bonnet calculations (see Santal\'o, \cite{san76}, p.112), to the following definition.
\begin{definition}\hspace{-2mm}\emph{\textbf{:}}
    In the graph $\calG$, the \textbf{Euler Entity $\chi$ of a face} $F$ is defined as the total (over all $F$'s face-circuits) of the turning-angle sums divided by $2\pi$,  \texttt{minus} the number of isolated vertices of valency zero in $F$'s interior. In the planar partitioning induced by the graph $\calG$, the\textbf{ Euler Entity $\chi$ of a cell} is the Euler Entity of the face in $\calG$ from which the cell is derived. The \textbf{Euler Entity $\chi$ of a cell-union} is the sum of the Euler Entities for the cells in the union.
\end{definition}
Thus the cell in Figure \ref{fig1} surrounding the $0$-valent vertex $P$ has $\chi=0$. Cells (a), (b), (c), (d), (f) and (h) have Euler Entities $-2, -1, 1, 1, -1$ and $0$ respectively. Cell-union (e) has $\chi=3$ provided the truncated right-most cell of this cell-union has no holes outside the window. The grey cell-union with one hole has $\chi = 1$.

Importantly, one must not interpret the face in Figure \ref{fig2} as having three holes. It has no holes; the only face-circuit in this face $F$ covers all of $F$'s boundary. Hence the domain has Euler Entity $\chi = 1$.

{\sc Sample results:} In this paper we provide natural generalisations for many of the geometric
and topological formulae given in the traditional theory cited in Section 1. The Euler Entity plays an important role. For example if, for the typical cell in a random planar partitioning (RPP) derived from a stationary ergodic random geometric graph, $\mu_\chi$ is the expected
Euler Entity, $\mu_E$ and $\mu_S$ are the cell's expected edge-count and expected side-count
whilst, for a typical vertex, $\theta$ is the expected valency and $\phi$ is the expected number of cell-side interiors containing the vertex, then we show
that
\begin{equation}\label{tag1}
\mu_E = \frac{2\theta\mu_\chi}{\theta-2}\quad\text{    and   }\quad \mu_S = \frac{2(\theta-\phi)\mu_\chi}{\theta-2},
\end{equation}
provided $\theta \neq 2$. We
further prove that $\theta=2$ if and only if $\mu_\chi = 0$. The formulae can also be adapted to studies of cell-unions.
\begin{remark}\hspace{-2mm}\emph{\textbf{:}}\label{thetais2}
     Consider a \tes\ comprising a lattice of regular hexagons with a vertex of valency $0$ placed at the centre of each hexagon, the whole structure being made stationary by randomising the planar origin within one hexagon. It provides an example of a generalised \tes\ having $\th=2$. Note that, because each cell has one hole, then $\mu_\chi = 0$. Another example with $\th=2$ and $\mu_\chi =0$ arises if each $0$-valent vertex in the example above is replaced by a short closed line-segment that does not hit any hexagon boundary. The ends of the line-segment produce two $1$-valent vertices.
\end{remark}

{\sc Method:} A key tool in our
proofs of these and other similar results, is the ergodic method described in \cite{cow78}-\cite{cow80} together with the following simple
generalisation of Euler's planar graph identity.

\begin{lemma}\hspace{-2mm}\emph{\textbf{:}}\label{genEuler}
    Consider a finite (not necessarily connected) planar graph $\bfG$  having $n$
nodes and $\ell$ links. We assume also that $\bfG$'s bounded faces has faces have a defined Euler Entity\footnote{$\bfG$ is finite --- not to be confused with the infinite graph $\calG$ discussed earlier. Unlike $\calG$, the finite graph $\bfG$ has an unbounded face.}.  Let $ \mathcal{X}$ be the sum of Euler Entities
over all bounded faces.  Then,
\begin{equation}\label{tag2}
n - \ell +   \mathcal{X} = 1.
\end{equation}
When the graph $\bfG$ is connected, all bounded faces have
$\chi=1$ and so then $  \mathcal{X}$ = the number of bounded faces
$f$. The familiar form of Euler's identity is $n-\ell+f=1$, often written as $n-\ell+f^* = 2$ where $f^* = f+1$, counting the unbounded face too.\hfill$\square$
\end{lemma}
The proof is by induction on
$\ell$, commencing with any case where $\bfG$ is connected.

\noindent
\section*{5. Motivating examples}\label{sec5}

\begin{figure}[h]
    \psfrag{c}{$\pi$}
    \begin{center}
        \includegraphics[width=75mm]{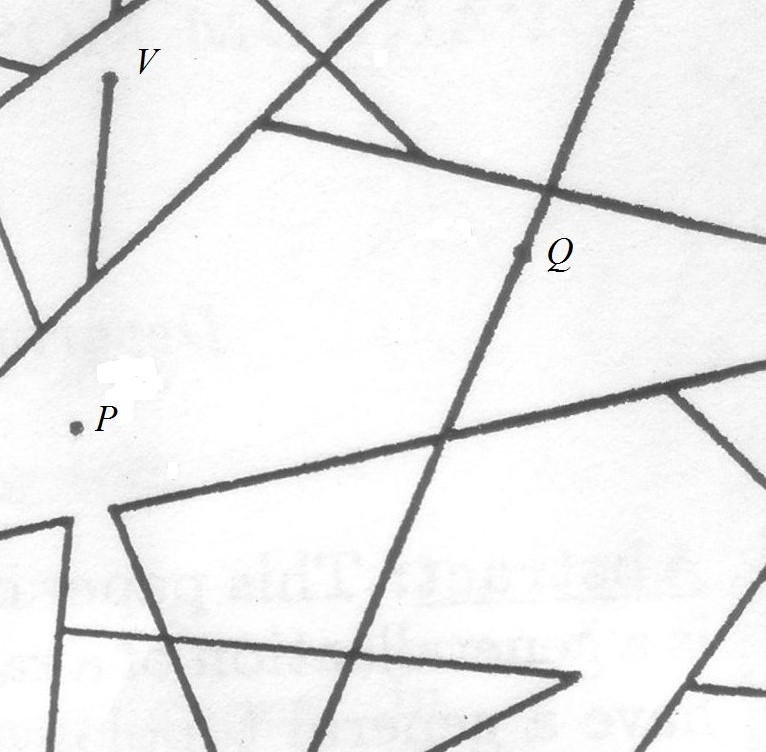}\hspace{3mm}
        \includegraphics[width=75mm]{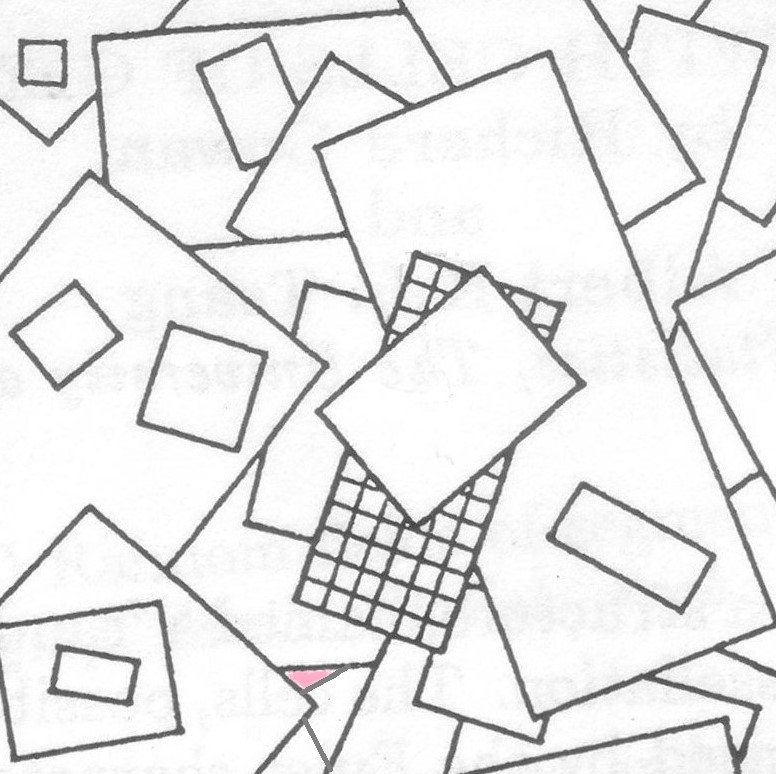}
        \\
        {\scriptsize (a)\hspace{7cm}(b)}
        \caption{\label{fig3} \scriptsize
(a) A \tes\ formed by a random line process with the random deletion of edge interiors.  Here $P$ is a $0$-valency vertex caused by the chance deletion of all its emanating edge-interiors.  (b) A realisation of a falling-leaf \tes\ with variation in the size and shape of the rectangular leaves.}
    \end{center}
\end{figure}

{\sc Random deletion of edge interiors:} Any \tes, even one of the traditional kind, may be altered by random edge
deletions, each edge interior  being deleted \ind ly with probability $q$.    Provided $p:= 1-q$ is not too small, the result of the deletions provides the frame of a generalised \tes,
including the possible creation of a 0-vertex (see $P$ in Figure
\ref{fig3}(a)). From this frame the open \tes-cells can be constructed. The figure, based on an initial stationary and isotropic Poisson line-process,
shows a $1$-vertex (see $V$), four $2$-vertices (one of these, $Q$, being of the collinear double-$\pi$ form) and numerous \pivs.

We have not yet investigated the critical value of $p$ for such line processes, below which the cells become unbounded. Readers will note that this structure is similar to structures studied in percolation theory where, in some cases, the critical value of $p$ is known. It will also be noted  that sometimes the percolation problem is cast as a random \emph{addition} of links to a stationary point process of nodes that initially has no links.

{\sc Falling leaf model:} When opaque leaves fall randomly on the plane \cite{cots94}, those falling later
will cover the others below.  Under assumptions of stationarity of
 leaves, the uncovered leaf-boundaries will
form a RPP.  The case of leaves congruent to a given simple polygon, where no leaf fits inside
another, was studied in \cite{cots94}\footnote{An isotropic assumption was also used in \cite{cots94}.}. With variation in size and shape of
leaves, however, cells of the \tes\ may have other cells wholly enclosed. So we have examples of the `holes' in cells (see Figure
\ref{fig3}(b) where the falling leaves are rectangles, assumed closed).

This model also provides a motivation for introducing \emph{cell-unions}.  The hatched domain in
Figure \ref{fig3}(b), comprises three disconnected cells. There is a clear nexus between these pieces; they all belong
to the same fallen leaf but have become disconnected by the position(s) of
a later leaf (or leaves). Perhaps one could declare that these cells be grouped as a cell-union. If so, $\chi = 3$ for the particular cell-union.

Another aspect of this falling-leaf \tes\ is the emergence of some closed cells, some open cells and some `neither open nor closed' cells. This occurs because if a closed leaf $L_1$ is first hit by a closed leaf $L_2$, without the boundary of $L_1$ being covered by $L_2$, then the visible part of $L_1$ is now $L_1\setminus L_2$ (and this is not open). In general, the visual part of a leaf whose boundary is partly covered, is \emph{neither open nor closed}. Yet in other situations, the visual part may be open or closed. A cell (shaded pink) belonging to a leaf whose boundary has been completely covered is an \emph{open} set. Recently-fallen closed leaves not yet hit by any later leaf are \emph{closed} sets. So, in order to conform to our theory, we need to focus on the new \tes\ frame at the time of observation and construct \emph{open} cells only from the frame (not a mix of topological types from the physical process of leaf-coverage). The falling leaf model also has many vertices of valency two and many \pivs.

\begin{figure}[h]
    \psfrag{c}{$\pi$}
    \begin{center}
        \includegraphics[width=70mm]{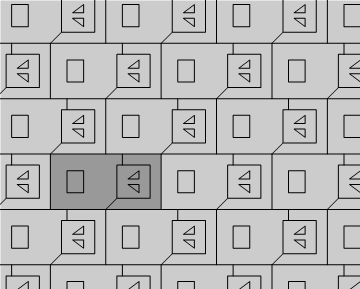}\hspace{5mm}
        \includegraphics[width=70mm]{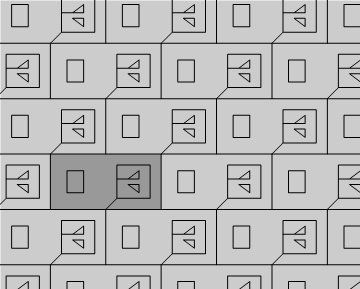}
        \\
        {\scriptsize (a)\hspace{7cm}(b)}
        \caption{\label{fig4}{\scriptsize  Two periodic tilings based, in each case, on the darkly shaded sub-collection of tiles (which occupy a rectangle of dimension $2\times 1)$ units. We assume the origin is uniformly distributed inside this rectangle, thereby making each tiling/\tes\ stationary.}}
    \end{center}
\end{figure}

{\sc Tilings beyond the theory of Gr\"{u}nbaum and Shephard:} In Figure \ref{fig4} are two tilings that have periodic repetition, via translation of a sub-collection of the tiles. These tilings violate the basic assumption N.1 of Gr\"{u}nbaum and Shephard \cite{grsh87} (which states that \emph{all tiles must be isomorphic to a closed disk}). Although the cells in Figure \ref{fig4}(a) might be considered \emph{closed}, doing so would lead to a violation of the assumption N.2 of these authors (\emph{the intersection of any two tiles is a connected set}); see page 121 of their book for these assumptions. With Figure \ref{fig4}(b), one of the cells cannot be considered \emph{closed} as this operation destroys a \tes\ edge.  We, of course, consider all cells in both figures as \emph{open} sets.

Our results bring some quantitative tools to these non-traditional tilings.  For the example in Figure \ref{fig4}(a), consider the six open cells in the darkly shaded rectangle viewed as follows:
\begin{itemize}
  \item a small rectangle (side-count $S= 4$, edge-count $E= 4$, Euler Entity $\calX =1$);
  \item a heptagon with the above rectangle as a hole ($S= 11$, $E=13$, $\calX = 0$);
  \item  two triangles (each has $S= 3$, $E=3$, $\calX = 1$);
  \item the rectangle which has the triangles as holes ($S= 10$, $E=11$, $\calX = -1$);
  \item the octagon ($S= 8$, $E=8$, $\calX = 1$).
\end{itemize}
So, by direct calculation, $\mu_S = (4 + 11 + 3 + 3 + 10 + 8)/6 = \tf{13}{2}$, $\mu_E = (4 + 13 + 3 + 3 + 11 + 8)/6 = 7$ and $\mu_\chi = (1 + 0 + 1 + 1 -1 + 1)/6 = \tf12$. Or, we might first observe that $18$ vertices are associated with the dark rectangle: the $15$ inside and the three left-most vertices on its boundary (we can't count all the boundary vertices as this would introduce much double counting across the \tes). From these $18$ vertices we calculate that $\th = \tf73$ and $\phi = \tf16$. Thus, from our new formulae in (\ref{tag1}),
\begin{equation*}
\mu_E = \frac{2\times \tf73 \times \tf12}{\tf73-2}= 7\quad\text{    and   }\quad \mu_S = \frac{2\times (\tf73 - \tf16)\times \tf12}{\tf73-2}=\frac{13}{2}.
\end{equation*}
Thus the directly calculated values are in accord with the formulae of our new theory shown in (\ref{tag1}). Analysis of Figure \ref{fig4}(b) is left as an exercise for the reader.

\smallskip
\noindent
\section*{ 6. Theoretical framework}

Our generalised theory  can follow \cite{cow80} to some extent, though initially we assume that no \emph{cell-unions} in the structure are of interest; thus we are only interested in cells.  As in Section 1, we let $(\Omega,
\calQ , \Pr )$ be a probability field, where $\Omega$ is the set of all allowed infinite
structures $\calG_\cup^+$ that are derived from an infinite geometric graph $\calG$. $ \calQ $ is a $\sigma$-algebra containing all the events
of interest to us.  The element $\omega \in \Omega$ is one realisation of $\calG_\cup^+$.
So $(\Omega, \calQ ,\Pr )$ is our random partition process.  For $x,\ t \in \zwei$ define the
translation operator $T^t: x \rightarrow x + t$.  Thus $T^t$, defined on
$\Omega$, translates any realisation $\omega$ by $t$.  As in Section 1, we assume that the
process is stationary and ergodic, via the assumption that $T^t$ (now
defined on $ \calQ$) is a measure-preserving and ergodic operator.

An important consequence of ergodicity is
 Wiener's ergodic theorem \cite{cow78, cow80}, which states that, if $X$
is any random variable derived from the RPP such that $\Bbb E|X| < \infty$ and $B_r$ is the closed
ball of radius $r$, centre 0 in $\zwei$, then for almost all $\omega$
\begin{equation}\label{tag3}
\lim_{r\rightarrow\infty}\frac1{\pi r^2}\int_{B_r}X(T^t\omega)dt = \Bbb E(X).
\end{equation}

Let $D$ be any compact convex reference domain in $\zwei$ unrelated to our RPP.  In $D$, the `entire' edges or cells are those
{\it wholly} in $D$; the other edges and other cells that intersect $D$ are called 'truncated'. An `edge-part' or `cell-part' in $D$ refers respectively to any entire or truncated edge or cell.
A cell 'centre' is any convenient reference point of the cell, for
example the cell's centroid. Define

\begin{align*}
N(D)  := \ & \text{the number of cell-parts in $D$,}\\
N'(D)  := \ & \text{the number of entire cells in $D$,}\\
n_{\text{cells}}(D)  := \ & \text{the number of cell centres in $D$,}\\
\ell(D)  := \ & \text{the total length of edge-segments in $D$,}\\
M(D)  := \ & \text{the number of edge-parts in $D$,}\\
M'(D)  := \ & \text{the number of entire edges in $D$,}\\
M^ \partial (D)  := \ & \text{the number of hittings of $\partial D$ by
edges,}\\
n_{\text{verts[k]}}(D)  := \ & \text{the number of vertices within $D$ of valency $k$,}\\
n_{\text{edges}}(D)  := \ & \text{the number of edge mid-points in $D$,}\\
n_{\text{edges}}^*(D)  := \ & \text{the number of edge ends in $D$,}\\
n_{\text{$\pi$-verts[k]}}(D)  := \ & \text{the number of $\pi$-vertices in $D$ of
valency $k$.}\end{align*}
Note that a common symbol, a subscripted $n$, is used for counts of `points' (where the points might be vertices, cell centroids, edge mid-points or edge ends). The subscript indicates the type of point.

Where
there is a need to emphasise the dependence on the realisation $\omega$,
we use the extended notation, $\ell(D,\omega)$ say.
It can be shown, by some elementary inequalities mostly given in \cite{cow80}, that
the assumptions $\Bbb EM(D) < \infty$ and $\Bbb
En_{\text{verts[0]}}(D) < \infty$ are sufficient to ensure that all these quantities have
finite expectation for bounded $D$.  We make these assumptions and also
assume that
$\Bbb EM(D) > 0$ if $|D|$, the Lebesgue measure of $D$, is positive. Note
that
\begin{equation}\label{tag4}
n_{\text{edges}}^*(D) = \sum k\ n_{\text{verts[k]}}(D).
\end{equation}
This summation, and all summations in the remainder of this paper, are for $k\geq 0$, unless otherwise marked.

Under stationarity, $\Bbb E\ell(\cdot), \Bbb E n_{\text{cells}}(\cdot),\Bbb En_{\text{edges}}(\cdot),\Bbb
En_{\text{verts[k]}}(\cdot)$ and $\Bbb
En_{\text{$\pi$-verts[k]}}(\cdot)$ are measures, proportional to Lebesgue measure.  So we
may introduce the finite constants $\alpha, \lam_{\text{cells}},\lam_{edges},\lam_{\text{verts}[k]}$ and
$\lam_{\text{$\pi$-verts}[k]},\ (k \geq 0)$ such that

\begin{alignat*}{2}
\Bbb E n_{\text{cells}}(D) =\ & \lam_{\text{cells}}|D|\qquad\quad & \Bbb E\ell(D) =\ & \alpha|D|\\
\Bbb En_{\text{edges}}(D) =\ & \lam_{\text{edges}}|D|\qquad\quad & \Bbb En_{\text{verts[k]}}(D) =\ & \lam_{\text{verts}[k]}|D|\\
\Bbb En_{\text{$\pi$-verts[k]}}(D) =\ & \lam_{\text{$\pi$-verts}[k]}|D|.
\end{alignat*}

These parameters, except $\alpha$, are the intensities of stationary
point processes in $\zwei$. We see from (\ref{tag4}) that $\Ex n_{\text{edges}}^*(\cdot)$ is
also a measure with $\Ex n_{\text{edges}}^*(D) = \sum k\ \lam_{\text{verts}[k]}|D|$. Note that $\alpha >
0$ and $\lam_{\text{$\pi$-verts}[0]} = \lam_{\text{$\pi$-verts}[1]} = 0$. We also assume that $\lam_{\text{cells}}$, $\lam_{\text{edges}}$,
$\lam_{\text{verts}}  :=  \sum\lam_{\text{verts}[k]}$ and $\lam_{\text{$\pi$-verts}}  :=  \sum\lam_{\text{$\pi$-verts}[k]}$ are positive.

\smallskip
\noindent
\section*{ 7. Ergodic theory}

If $D$ is taken to be the ball $B_r$, ergodic arguments can now be applied to
show that, for example, $\ell(B_r)/\pi r^2 \rightarrow \alpha$ almost
surely as $r\rightarrow\infty$.  To understand the detailed use of (\ref{tag3}),
take $B_r$ and $B_y,\ y \ll r$, and consider a random variable,
associated with $B_y$.  Then consider this random variable for the
translated disk $T^{-t}B_y$ and integrate over all $t$, firstly within
$B_{r-y}$ and then within $B_{r+y}$.  For example, consider the two integrals
$$
I_1 = \int_{B_{r-y}}n_{\text{verts[k]}}(T^{-t}B_y,\omega)dt,\qquad
I_2 = \int_{B_{r+y}}n_{\text{verts[k]}}(T^{-t}B_y,\omega)dt.$$
If we now draw a circle of radius $y$ around each $k$-vertex and consider
the sum of these circular areas (including any parts which may extend
beyond $B_r$) then this sum, which is obviously equal to $\pi y^2
n_{\text{verts[k]}}(B_r,\omega)$, is bounded below by $I_1$ and above by $I_2$.  Noting
that $n_{\text{verts[k]}}(T^{-t}B_y,\omega) = n_{\text{verts[k]}}(B_y,T^t\omega)$, we have that
$$
\int_{B_{r-y}}n_{\text{verts[k]}}(B_y,T^t\omega)dt \leqq \pi y^2n_{\text{verts[k]}}(B_r,\omega)
\leqq \int_{B_{r+y}}n_{\text{verts[k]}}(B_y,T^t\omega)dt.$$

Therefore
$$
\frac{\pi(r-y)^2}{\pi r^2}
\int_{B_{r-y}}\hspace{-4mm}\frac{n_{\text{verts[k]}}(B_y,T^t\omega)}{\pi(r-y)^2}dt \leqq \frac{\pi
y^2 n_{\text{verts[k]}}(B_r,\omega)}{\pi r^2}
\leqq \frac{\pi(r+y)^2}{\pi r^2} \int_{B_{r+y}}\hspace{-4mm}\frac{
n_{\text{verts[k]}}(B_y,T^t\omega)}{\pi(r+y)^2}dt.$$
The left and right sides of this inequality converge with probability one
to the same quantity, $\Ex n_{\text{verts[k]}}(B_y)$, by applying the Wiener
ergodic theorem (\ref{tag3}).  Therefore
\begin{equation}\label{tag5}
    \frac{n_{\text{verts[k]}}(B_r,\omega)}{\pi r^2} \overset{\text{a.s.}}{\longrightarrow}
\frac{\Ex n_{\text{verts[k]}}(B_y)}{\pi y^2} = \lam_{\text{verts}[k]}.
\end{equation}

With an almost identical argument one can show that the other counting
variates associated with point processes converge almost surely.
\begin{align}
\frac{n_{\text{$\pi$-verts[k]}}(B_r,\omega)}{\pi r^2} \overset{\text{a.s.}}{\longrightarrow}
\frac{\Ex n_{\text{$\pi$-verts[k]}}(B_y)}{\pi y^2} =\ & \lam_{\text{$\pi$-verts}[k]},\label{tag6}\\
\frac{n_{\text{edges}}B_r,\omega)}{\pi r^2} \overset{\text{a.s.}}{\longrightarrow}
\frac{\Ex n_{\text{edges}}B_y)}{\pi y^2} =\ & \lam_{\text{edges}},\label{tag7}\\
\frac{n_{\text{cells}}(B_r,\omega)}{\pi r^2} \overset{\text{a.s.}}{\longrightarrow}
\frac{\Ex n_{\text{cells}}(B_y)}{\pi y^2} =\ & \lam_{\text{cells}},\label{tag8}
\end{align}
whilst from (\ref{tag4}) and (\ref{tag5})
\begin{align}
\frac{n_{\text{edges}}^*B_r,\omega)}{\pi r^2} \overset{\text{a.s.}}{\longrightarrow}
\frac{\Ex n_{\text{edges}}^*B_y)}{\pi y^2} = \sum k\ \lam_{\text{verts}[k]}.\label{tag9}
\end{align}

If the quantity $\ell(B_y,\omega)$ is integrated in the same manner we
obtain the inequality
$$
\int_{B_{r-y}}\ell(B_y,T^t\omega)dt \leq \pi y^2\ell(B_r,\omega) \leq \int_{B_{r+y}}\ell(B_y,T^t\omega)dt.$$
Here, the middle expression involves a small calculation.  Around each
segment in $B_r$, construct a sausage-shaped domain of points within
distance $y$ of that segment.  As the centre of a disk $B_y$ moves over
all positions within the domain, it can easily be shown that an
integration of the segment length within $B_y$ yields $\pi y^2$ times the
segment length.  Adding over all `sausage' domains yields $\pi
y^2\ell(B_r,\omega)$.  Dividing by $\pi r^2$, taking limits and applying
the Wiener ergodic theorem proves that
\begin{equation}
    \frac{\ell(B_r,\omega)}{\pi r^2}\overset{\text{a.s.}}{\longrightarrow}
\alpha.\label{tag10}
\end{equation}

Next we consider the integrals $I_1$ and $I_2$ (say) of
$M(B_y,T^t\omega)$, or equivalently $M(T^{-t}B_y,\omega)$, as $t$ moves
over $B_{r-y}$ and $B_{r+y}$ respectively.  These integrals provide lower
and upper bounds for the sum of areas for all of the `sausage' domains.
This sum is easily seen to be $2y\,\ell(B_r,\omega)+\pi y^2M(B_r,\omega)$.  Thus
$$
I_1 := \int_{B_{r-y}}M(B_y,T^t\omega)dt \leq 2y\,\ell(B_r,\omega) + \pi
y^2M(B_r,\omega) \leq \int_{B_{r+y}}M(B_y,T^t\omega)dt =: I_2.$$
Dividing by $\pi r^2$, taking limits and using (\ref{tag3}) and (\ref{tag10}), shows that
\begin{equation}\label{tag11}
\frac{M(B_r,\omega)}{\pi r^2}\overset{\text{a.s.}}{\longrightarrow}
\frac{\Ex M(B_y)}{\pi y^2} - \frac{2\alpha}{\pi y}.
\end{equation}
In adding the areas of the `sausage' domains, we included the
semi-circular parts which extend beyond $B_r$ when an edge hits $\partial
B_r$.  If these semi-circular areas are not counted, we find that
\begin{equation}
I_1 \leq 2y\,\ell(B_r,\omega) + \frac12 \pi y^2\sum k\ n_{\text{verts[k]}}(B_r,\omega) \leq
I_2.\label{tag12}
\end{equation}
This is a precise way of saying that, since each edge has two ends,
$2M(B_r) = \sum k\ n_{\text{verts[k]}}(B_r)$ {\it except for} the boundary effects of
$\partial B_r$.  Thus, the middle expression in (\ref{tag12}), when divided by $\pi
r^2$ converges almost surely to $\Ex M(B_r)$.  We already know however from
(\ref{tag5}) and (\ref{tag10}), that it converges to $2y\alpha + \frac12\pi y^2\sum k\ \lam_{\text{verts}[k]}$,
so
\begin{equation}
\Ex M(B_y) = 2y\alpha + \frac12\pi y^2\sum k\ \lam_{\text{verts}[k]}.\label{tag13}
\end{equation}
From (\ref{tag11}) and (\ref{tag13}), therefore,
\begin{equation}
\frac{M(B_r)}{\pi r^2} \overset{\text{a.s.}}{\longrightarrow} \frac12 \sum
k\ \lam_{\text{verts}[k]}.\label{tag14}
\end{equation}

\noindent
\section*{ 8. Sampling the typical vertex or typical edge}

The mean valency of a `typical' vertex of the RPP, denoted by $\theta$, is defined as
the limit of the total valency of vertices within $B_r$ divided by the
number of vertices in $B_r$, as $r \rightarrow \infty$, whenever this
almost-sure limit exists and yields a constant.  We have established
existence because, using (\ref{tag5}),
\begin{align}
\theta  := \ & \lim_{r\rightarrow \infty}\frac{\sum
k\ n_{\text{verts[k]}}(B_r,\omega)}{\sum n_{\text{verts[k]}}(B_r,\omega)}\notag\\
=\ & \frac{\sum k\ \lam_{\text{verts}[k]}}{\sum\lam_{\text{verts}[k]}} = \frac{\sum k\ \lam_{\text{verts}[k]}}{\lam_{\text{verts}}}\label{tag15}
\end{align}
where $\lam_{\text{verts}}$ is the intensity of the point process of {\it all} vertices.

In a similar fashion, the mean length of a `typical' edge of the RPP, denoted by
$\nu$, is defined as the limit of total segment length within $B_r$,
divided by the number of edge-segments in $B_r$.  Thus, from (\ref{tag10}) and (\ref{tag14}),

\begin{align}
\nu  := \ & \lim_{r\rightarrow\infty}\frac{\ell(B_r,\omega)}{M(B_r,\omega)}\notag\\
=\ & \frac{2\alpha}{\sum k\ \lam_{\text{verts}[k]}} = \frac{2\alpha}{\lam_{\text{verts}}\theta}.\label{tag16}
\end{align}

\noindent
\section*{ 9. Edges hitting the boundary}

The edges of our process can be viewed as a stationary line-segment process
(LSP) of a general kind. Following \cite{cow79}, we have that for $y>0$
$$
M(B_r)-M'(B_r) \leq [\ell(B_{r+y})-\ell(B_r)]/y+n_{\text{edges}}^*(B_{r+y})-n_{\text{edges}}^*(B_r).$$

Thus from this inequality, combined with (\ref{tag9}) and (\ref{tag10}), we see that the
normalised number of edges hitting the boundary is almost surely
asymptotically negligible as $r\ra\i$, that is,
\begin{equation}
    \frac{M(B_r)-M'(B_r)}{\pi r^2} \overset{\text{a.s.}} \longrightarrow
0.\label{tag17}
\end{equation}

Since $M' \leq n_{\text{edges}} \leq M$, we see that $n_{\text{edges}}(B_r)/\pi r^2$ and $M(B_r)/\pi r^2$
converge almost surely to the same limit, namely that given in (\ref{tag14}). Thus
from (\ref{tag7}), (\ref{tag14}) and (\ref{tag15}), $2\lam_{\text{edges}} = \lam_{\text{verts}}\theta$, a result which, combined
with (\ref{tag16}), permits many rearrangements, for example,
\begin{equation}
    \lam_{\text{edges}} = \frac{\alpha}{\nu}.\label{tag18}
\end{equation}

Note that $M-M' = M^{\partial}_{1}+M^{\partial}_{2}$ where
$M^{\partial}_{i}$ is
defined as the number of edges which cut the boundary $i$ times. Thus
$M^{\partial}_{i}(B_r)/\pi r^2 \overset{\text{a.s.}}{\longrightarrow} 0$
and since $M^{\partial} = M^{\partial}_{1}+2M^{\partial}_{2}$, we see also
that $M^{\partial}(B_r)/\pi r^2 \overset{\text{a.s.}}{\longrightarrow} 0$.

\noindent
\section*{ 10. Sampling the typical cell}

We now focus attention on the random finite graph $  \mathbf{G_r(\om)}$ whose nodes and links are as follows.
\begin{itemize}
  \item The nodes are all vertices of the RPP which lie in $B_r$, together with all points where the boundary $\partial B_r$ intersects an edge of the RPP.
  \item The links are all edge-parts in $B_r$ together with all the circular arcs which make up $\partial B_r$.
\end{itemize}
Here $\om$ is the infinite graph realised randomly. Any cell-part formed within $B_r$ has an area $A$ and perimeter $L$
(perhaps involving part of $\partial B_r$). The sums of area and perimeter over all regions are denoted by $\mathcal{A}(B_r)$ and
 $\mathcal{L}(B_r)$ and it is clear that
\begin{equation*}
  \mathcal{A}(B_r) =\  \pi r^2\qquad\text{and}\qquad
 \mathcal{L}(B_r) =\  2\ell(B_r) + 2\pi r.
\end{equation*}

\begin{figure}[h]
    \psfrag{c}{$\pi$}
    \begin{center}
        \includegraphics[width=90mm]{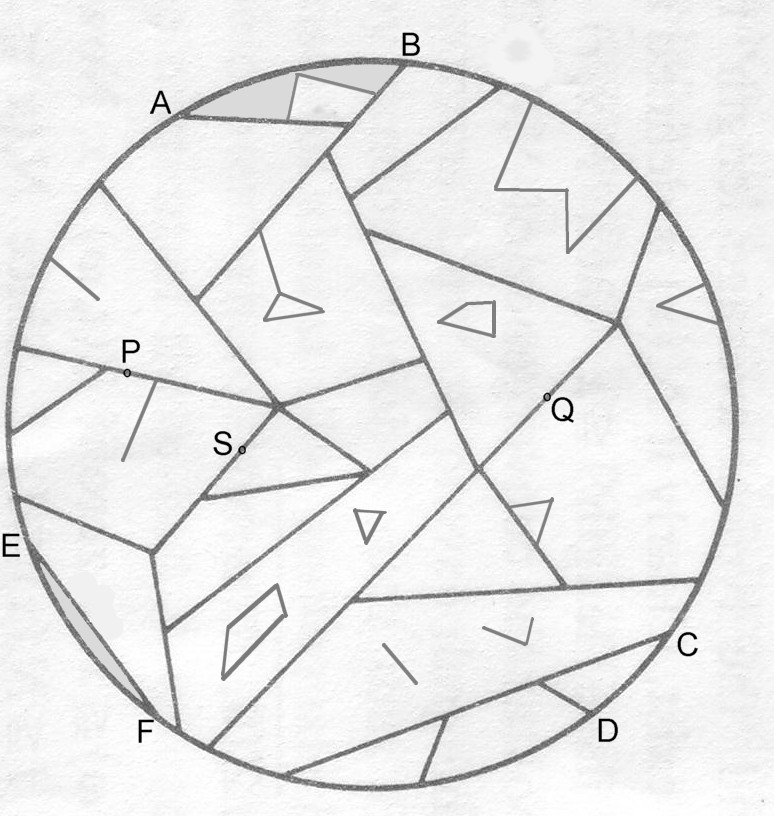}
        \\
        \caption{\label{fig5}{\scriptsize  The finite graph $\bfG_r(\om)$ comprises the planar partitioning within $B_r$ together with $\partial B_r$. All double-$\pi$ vertices, those labelled $\mathsf P, \mathsf Q$ and $\mathsf S$, are marked. For the entire cells, the Euler Entity, the edge-count, side-count, etcetera, are defined by our earlier definitions. For a truncated cell, which typically has an arc from $\partial B_r$ on its boundary, the text describes one how these quantities are calculated. For example, the cell-parts whose arcs are $\mathsf A\mathsf B, \mathsf C\mathsf D$ or $\mathsf E\mathsf F$ each have Euler Entity $1$.
        }}
    \end{center}
\end{figure}

The vertex-count $V$, edge-count $E$, side-count $S$ and corner-count $C$ of a cell-part are defined in Definition \ref{defn_edgecount} when the cell-part is an entire cell. Those definitions apply immediately to the truncated cells by simply treating the arcs as links. For example, in Figure \ref{fig5} the side-count of the two shaded cell-parts  are $5$ and $2$, treating the arcs $\mathsf A \mathsf B$ and $\mathsf E\mathsf F$ as sides. One can also extend the definition of the Euler Entity to any truncated cell having circular arcs on its boundary; we add an `arc turning angle' which, in the spirit of Gauss-Bonnet, is the total angle turned by the walker when he traverses the arc. It can be expressed as an integral of curvature over the arc even for open sets (see Santal\'o, \cite{san76}, formula 7.16).

\begin{remark}\colon
For most of the truncated cells, the Euler Entity can be calculated correctly by first replacing any arc, say $\mathsf C\mathsf D$, by a \lseg\ having the same end points. But this device doesn't work when the truncated cell is like one of the two shaded sets. In these cases, however, one can replace the arc with a polygonal chain --- then calculate the Euler Entity using our usual definition for simple polygons. The use of an \emph{arc turning angle} is the simplest approach, we think.
\end{remark}

Sums of the entities vertex-count $V$, edge-count $E$, side-count $S$ and corner-count $C$ over all cell-parts in $B_r$ are denoted by $   \mathcal{V}(B_r), \mathcal{E}(B_r), \mathcal{S}(B_r)$ and $\mathcal{C}(B_r)$, using script letters.
Clearly,

\begin{align}
 \mathcal{E}(B_r) =\ & 4M(B_r) - \sum k\ n_{\text{verts[k]}}(B_r)\notag\\
  \mathcal{C}(B_r) =\ &  \mathcal{E}(B_r)-2n_{\text{$\pi$-verts[2]}}(B_r) - \sum_{k\geq 3}n_{\text{$\pi$-verts[k]}}(B_r)\notag\\
  \mathcal{V}(B_r) =\ &  \mathcal{} E(B_r) + n_{\text{verts[0]}}(B_r)\label{tag19}\\
  \mathcal{S}(B_r) =\ &   \mathcal{C}(B_r).\notag
\end{align}

Moreover the sum $  \mathcal{X}(B_r)$ of Euler Entities over all cell-parts
can be found from Lemma \ref{genEuler} using the finite graph $  \mathbf{G}_r(\om)$.  Thus $n =
\sum n_{\text{verts[k]}}(B_r) + M^{\partial}(B_r),\ \ell = M(B_r) + M^{\partial}(B_r)$ and
$  \mathcal{X} =   \mathcal{X}(B_r)$, so (\ref{tag2}) becomes
\begin{equation}
     \mathcal{X}(B_r) = M(B_r) - \sum n_{\text{verts[k]}}(B_r) + 1.\label{tag20}
\end{equation}

We define the mean of a cell-part feature, for example of area $A$, as the
limit as $r\rightarrow \infty$ of an appropriate cell-part sum, $ \mathcal{}
A(B_r)$ say, divided by the number $N(B_r)$ of cells, {\it if} this
ratio converges to a constant almost surely.  We denote such mean values
by $\mu$, suitably subscripted.  Thus

\begin{align}
\mu_A  := \ & \lim_{r\rightarrow\infty}\frac{  \mathcal{A}(B_r)}{N(B_r)}\notag\\
\mu_L  := \ & \lim_{r\rightarrow\infty}\frac{ \mathcal{L} (B_r)}{N(B_r)}\notag\\
\mu_\chi  := \ & \lim_{r\rightarrow\infty}\frac{  \mathcal{X}(B_r)}{N(B_r)}\notag\\
\mu_V  := \ & \lim_{r\rightarrow\infty}\frac{  \mathcal{V}(B_r)}{N(B_r)}\label{tag21}\\
\mu_E  := \ & \lim_{r\rightarrow\infty}\frac{ \mathcal{E} (B_r)}{N(B_r)}\notag\\
\mu_S = \mu_C  := \ & \lim_{r\rightarrow\infty}\frac{  \mathcal{C}(B_r)}{N(B_r)}\notag
\end{align}
whenever the appropriate limit exists almost surely.  In \cite{cow80}, these
entities were defined as the limit of the {\it expected} ratio, for
example $\lim \Ex(  \mathcal{A}(B_r)/N(B_r))$, when that limit exists.  The
definition that we adopt in (\ref{tag21}), and earlier in (\ref{tag15}) and (\ref{tag16}), avoids
certain technicalities and is, in all examples that we have experienced,
equivalent to that in \cite{cow80}.

Thus we can say, from (\ref{tag5}), (\ref{tag6}), (\ref{tag10}) and (\ref{tag14}) that
\begin{align}
\frac{  \mathcal{A}(B_r)}{\pi r^2}\ & \overset{\text{a.s.}}\longrightarrow
1\notag\\
\frac{  \mathcal{X}(B_r)}{\pi r^2}\ & \overset{\text{a.s.}}\longrightarrow
\frac12\sum^\infty_{k=0}(k-2)\lam_{\text{verts}[k]} = \frac{\lam_{\text{verts}}(\theta-2)}2
\notag\\
\frac{ \mathcal{L} (B_r)}{\pi r^2}\ & \overset{\text{a.s.}}\longrightarrow
2\alpha\notag\\
\frac{  \mathcal{V}(B_r)}{\pi r^2}\ & \overset{\text{a.s.}}\longrightarrow
\lam_{\text{verts}[0]} + \sum^\infty_{k=1}k\ \lam_{\text{verts}[k]} = \lam_{\text{verts}}(\theta+\xi)\label{tag22}
\\
\frac{ \mathcal{E} (B_r)}{\pi r^2}\ & \overset{\text{a.s.}}\longrightarrow
\sum^\infty_{k=1}k\ \lam_{\text{verts}[k]} = \lam_{\text{verts}}\theta \notag\\
\frac{  \mathcal{S}(B_r)}{\pi r^2}\ & \overset{\text{a.s.}}\longrightarrow
\sum^\infty_{k=1}k\ \lam_{\text{verts}[k]} - \sum^\infty_{k=3}\lam_{\text{$\pi$-verts}[k]} - 2\lam_{\text{$\pi$-verts}[2]} =
\lam_{\text{verts}}(\theta-\phi).\notag
\end{align}
Here $\phi$ is the mean number of angles equal to $\pi$ at a `typical'
vertex.  Formally, $\phi$ is the almost-sure limit
$$
\lim_{r\rightarrow\infty}\frac{\sum\limits^\infty_{k=3}n_{\text{$\pi$-verts[k]}}(B_r) +
2n_{\text{$\pi$-verts[2]}}(B_r)}{\sum n_{\text{verts[k]}}(B_r)} = \frac{ \lam_{\text{$\pi$-verts}} + \lam_{\text{$\pi$-verts}[2]}}{\lam_{\text{verts}}}.$$
In addition, $\xi$ is the proportion of vertices which are of valency zero.

These results show that the numerators in (\ref{tag21}) converge almost surely to constants,
when normalised by $\pi r^2$.  We now show that $N(B_r)/\pi r^2$ converges
likewise,
 thereby establishing the conditions for
the mean features for typical cells to be finite.

\noindent
\section*{ 11. Asymptotics of $N(B_r)$}

Firstly we need to show that $N(B_r) - N'(B_r)$, the number of truncated
cells, becomes asymptotically negligible relative to $\pi r^2$ as
$r\rightarrow\infty$.  It is a trivial fact that the number of
truncated cells is bounded above by $M^\partial(B_r)$.  For an ergodic \lseg\ process it is
shown in \cite{cow79} that, provided the expected number of these line-segments
hitting a bounded domain $D$ is finite, the number of crossing points of
edges with $\partial B_r$, when normalised by $\pi r^2$, tends almost
surely to zero.  Thus in our theory, $M^\partial(B_r)/\pi r^2
\overset{\text{a.s.}}\longrightarrow 0$ since we already have the
regularity condition $\Ex M(D) < \infty$.
Thus
\begin{equation}
    \frac{N(B_r) - N'(B_r)}{\pi r^2}\overset{\text{a.s.}}\longrightarrow 0.\label{tag23}
\end{equation}
Let $A_t$ be the area of the cell which covers $t\in \zwei$ with $A_t$ defined
as zero if $t$ lies on an edge of the RPP.  Stationarity implies that the
distribution of $A_t$ is independent of $t$.  Now, following \cite{cow80}, consider
the integral
$$
I = \int_{B_r}\frac{dt}{A_t(\omega)} = \int_{B_r}\frac{dt}{A_0(T^{-t}\omega)}.$$
This integral is approximately equal to $N(B_r,\omega)$.  Precisely
$N'(B_r,\omega) \leq I \leq N(B_r,\omega)$.  Since $N(B_r) \leq M(B_r)+1$,
$EN(B_r) < \infty$, so $E(I)$ is finite.  Thus from Fubini's theorem and
homogeneity, $E(1/A_0)$ is finite.  Wiener's theorem can thus be employed
to show that $I/\pi r^2 \overset{\text{a.s.}}{\longrightarrow} E(1/A_0)$.
 Rewriting the inequality as $I-(N-N') \leq N'\leq I$ and noting (\ref{tag23}), we
establish that $N'/\pi r^2\overset{\text{a.s.}}{\longrightarrow}
E(1/A_0)$ which, from (\ref{tag23}) too, implies that
\begin{equation}
    \frac{N(B_r)}{\pi r^2} \overset{\text{a.s.}}{\longrightarrow} E(1/A_0).\label{tag24}
\end{equation}

\smallskip
\noindent
\section*{ 12. Cellular mean values}

We have established that the denominators in (\ref{tag16}), normalised by
$\pi r^2$, converge almost surely to a constant.  Thus in conducting the
`ergodic experiment' to sample the typical cell, we have proved the
finiteness of mean values $\mu_A,\ \mu_L,\ \cdots$, defined in (\ref{tag21}).  In particular
\begin{equation}
    \mu_A = \lim_{r\rightarrow\infty}\frac{  \mathcal{A}(B_r)/\pi r^2}{N(B_r)/\pi
r^2} = \frac1{E(1/A_0)}.\label{tag25}
\end{equation}

Thus the unfamiliar entity $E(1/A_0)$, which appears in (\ref{tag24}), has a
convenient evaluation in terms of the mean cell area, namely
\begin{equation}
    E(1/A_0) = 1/\mu_A.\label{tag26}
\end{equation}
Further results which follow directly from (\ref{tag21}), (\ref{tag22}) and (\ref{tag24}) are, using
(\ref{tag26}), (\ref{tag16}) and (\ref{tag18}),

\begin{align}
\mu_L =\ & 2\alpha\mu_A = \lam_{\text{verts}}\theta\nu\mu_A = 2\lam_{\text{edges}}\nu\mu_A\notag\\
\mu_\chi =\ & \frac{\lam_{\text{verts}}(\theta-2)\mu_A}2 = (\lam_{\text{edges}}-\lam_{\text{verts}})\mu_A\label{tag27}
\end{align}

\begin{align}
\mu_E =\ & \lam_{\text{verts}}\theta\mu_A = 2\lam_{\text{edges}}\mu_A\notag\\
\mu_S = \mu_C =\ & \lam_{\text{verts}}(\theta-\phi)\mu_A\notag\\
\mu_V =\ & \lam_{\text{verts}}(\theta+\xi)\mu_A\label{tag28}
\end{align}
Clearly these formulae permit a large number of rearrangements including
the interesting topological-linkage formulae promised in (\ref{tag1}).
\begin{align}
\mu_E =\ & \frac{2\theta\mu_\chi}{\theta-2}\label{tag29}\\
\mu_S = \mu_C =\ & \frac{2(\theta-\phi)\mu_\chi}{\theta-2}\label{tag30}\\
\mu_V =\ & \frac{2(\theta+\xi)\mu_\chi}{\theta-2}\label{tag31}
\end{align}
which hold when $\theta \neq 2$.

Using (\ref{tag27}), it is readily proved that $\theta=2$ if and only if
$\mu_\chi = 0$, since $\lam_{\text{verts}}>0$ and $\mu_A > 0$.  Formulae (\ref{tag29})--(\ref{tag31}) are interesting
generalisations of formulae for the ergodic convex-celled \tes\ treated
in \cite{cow78}, \cite{cow80} and \cite{como88}, where $\mu_\chi=1$.
\begin{align}
\mu_E =\ & \mu_V = \frac{2\theta}{\theta-2}\label{tag32}\\
\mu_S =\ & \mu_C = \frac{2(\theta-\phi)}{\theta-2}.\label{tag33}
\end{align}

\smallskip
\noindent
\section*{ 13. Point processes}

We have already seen that there is a point process of vertices, intensity
$\lam_{\text{verts}}$, and a point process of edge mid-points, intensity $\lam_{\text{edges}}$.  Cells
can be given a reference point, for example the centroid, and these
reference locations form a point process, whose intensity we denote by
$\lam_{\text{cells}}$.  With all point processes, the count of points within $B_r$,
divided by $\pi r^2$, has an almost-sure limit equal to the intensity as
$r\rightarrow\infty$, under the ergodicity assumption.  This can be shown
using the methodology leading to (\ref{tag5})--(\ref{tag8}).

The choice of reference point is somewhat arbitrary and for our current
purpose it is convenient to choose a reference point which always lies in
the topological closure of the cell.  Centroids may not, so we choose the
mid-point of the longest cell-side.  Let $n_{\text{cells}}(B_r)$ be the number of cell
reference points in $B_r$.  On the one hand $n_{\text{cells}}(B_r)/\pi r^2
\overset{\text{a.s.}}{\longrightarrow} \lam_{\text{cells}}$.  On the other hand, $N'(B_r)
\leq n_{\text{cells}}(B_r) \leq N(B_r)$, so from (\ref{tag22}), (\ref{tag24}) and (\ref{tag26}), $n_{\text{cells}}(B_r)/\pi r^2
\overset{\text{a.s.}}{\longrightarrow} 1/\mu_A$.  Therefore
$$
\lam_{\text{cells}} = \frac1{\mu_A}.$$
Substitution for $\mu_A$ in (\ref{tag27}) yields
$$
\lam_{\text{edges}} = \lam_{\text{cells}}\mu_\chi + \lam_{\text{verts}}.$$
This result generalises the classical formula linking the three
point-process intensities.  Classically, $\chi = 1$ for all cells and so
$\lam_{\text{edges}} = \lam_{\text{cells}} + \lam_{\text{verts}}$ (as first shown in \cite{mec80} for tessellations
containing only convex cells).

\smallskip
\noindent
\section*{ 14. Ignoring vertices of valency 2}

In this section we generalise a result, first mentioned by Miles \cite{mil88},
involving a special type of vertex counting.  Let $\theta^\ast$ be the
expectation of a typical vertex's valency {\it conditional} upon the valency
{\it not} being equal to two.  Let $\mu_{V^\ast}$ be the mean number of
vertices for a typical cell ignoring vertices of valency 2.  For
tessellations where $\lam_{\text{verts}[0]} = \lam_{\text{verts}[1]} = 0$ and where each cell has $\chi
= 1$, it is argued by Miles that
\begin{equation}
    \mu_{V^\ast} = \frac{2\theta^\ast}{\theta^\ast-2}.\label{tag34}
\end{equation}

Within our more generalised RPP structure, we formally define

\begin{align*}
\theta^\ast =\ & \lim_{r\rightarrow\infty}\frac{\sum k\ n_{\text{verts[k]}} -2 n_{\text{verts[2]}}}{\sum n_{\text{verts[k]}}
- n_{\text{verts[2]}}}\\
=\ & \frac{\lam_{\text{verts}}\theta-2\lam_{\text{verts}[2]}}{\lam_{\text{verts}}-\lam_{\text{verts}[2]}}\\
\mu_{V^\ast} =\ & \lim_{r\rightarrow\infty}\frac{\nu(B_r) -2 n_{\text{verts[2]}}(B_r)}{N(B_r)}\\
=\ & \mu_V - 2\lam_{\text{verts}[2]}\mu_A\\
=\ & \frac{2[\lam_{\text{verts}}(\theta+\xi)-2\lam_{\text{verts}[2]}]\mu_\chi}{\lam_{\text{verts}}(\theta-2)}
\end{align*}
using (\ref{tag27}) and (\ref{tag28}).  Some rearrangement yields
\begin{equation}\label{tag35}
\mu_{V^\ast} = \frac{2(\theta^\ast+\xi^\ast)}{\theta^\ast-2}\mu_\chi
\end{equation}
where $\xi^\ast = \lam_{\text{verts}[0]}/(\lam_{\text{verts}}-\lam_{\text{verts}[2]})$, the proportion of zero-valency
vertices when $2$-valent vertices are ignored.  This formula (\ref{tag35}) is a
precise analogy of (\ref{tag31}), and generalises (\ref{tag34}).  Miles does not comment on
the mean number of corners when valency-2 vertices are ignored.
Let $\mu_{C^\ast}$ be this conditional mean.
\begin{align*}
\mu_{C^\ast} =\ & \lim_{r\rightarrow\infty} \frac{  \mathcal{C}(B_r) +
2  n_{\text{$\pi$-verts}[2]}(B_r) - 2n_{\text{verts[2]}}(B_r)}{N(B_r)}\\
=\ & \mu_C + (2\lam_{\text{$\pi$-verts}[2]}-2\lam_{\text{verts}[2]})\mu_A\\
=\ & \frac{2[\lam_{\text{verts}}(\theta-\phi)+2(\lam_{\text{$\pi$-verts}[2]}-\lam_{\text{verts}[2]})]\mu_\chi}{\lam_{\text{verts}}(\theta-2)}.
\end{align*}
Rearrangement yields a formula analogous to one of the results in (\ref{tag30}),
$$
\mu_{C^\ast} = \frac{2(\theta^\ast-\phi^\ast)}{\theta^\ast-2}\mu_\chi$$
where $\phi^\ast$ is the conditional mean number of angles equal to $\pi$
in a typical vertex, namely $( \lam_{\text{$\pi$-verts}}-\lam_{\text{$\pi$-verts}[2]})/(\lam_{\text{verts}}-\lam_{\text{verts}[2]})$ or $(\lam_{\text{verts}}\phi-2\lam_{\text{$\pi$-verts}[2]})/(\lam_{\text{verts}}-\lam_{\text{verts}[2]})$.

Nothing of interest happens when vertices of other valencies are ignored;
here vertices of valency two have a special status.

\section*{15. Extension of the ideas to cell-unions}

We have found, partly through experimentation, that the formulae in our theory can be applied to cell-unions (instead of to cells alone).  For example, the formulae in (\ref{tag1}) are valid if $\mu_E, \mu_S$ and $\mu_\chi$ are redefined as the expected edge-count of the cell-union, the expected side-count of the cell-union and the expected Euler entity of the cell-union. To reinforce this, we use bold fonts, $\pmb\mu_E, \pmb\mu_S$ and $\pmb\mu_\chi$, when calculating for expected values of cell-unions.

By some well-defined unionisation rule, cells are grouped --- the union is taken of those in each group. The rule is such that all groups contain a finite expected-number of cells.  Some cells may not be involved in a union; they are then in a `group of size one'.

We do not present any formal theory here, as there are many situations and many ways that unions of cells might be made. So formal arguments that embrace all possibilities are left to later publications.

In this paper we merely demonstrate the ideas, using the examples in Figure \ref{fig4}. Two unionisation rules are given.
\begin{itemize}
    \item A: The six cells in the dark region are grouped as follows: the two triangles form a group (whose union has $\chi = 2, E = 6$ and $S = 6$) and the two cells with rectangular outer-boundaries form a group (whose union has $\chi = 0, E = 15$ and $S = 14$). So there are four cell-unions in the dark region. The same grouping is applied periodically to all copies of the dark region. Using cell-unions rather than cells, $\pmb\mu_E = \tfrac{42}{4},\  \pmb\mu_S=\tfrac{39}{4}$ and $\pmb\mu_\chi = \frac34$.
    \item B: The six cells in the dark region are grouped according to their Euler Entity. So four cells with $\chi = 1$ make up one group (whose cell-union has $\chi = 4, E = 18$ and $S = 18$). Two other groups have just one cell, the heptagon with a rectangular hole ($\chi=0, E = 13, S = 11$) and the rectangle with two triangular holes ($\chi = -1, E = 11, S = 10$). So there are three cell-unions. Also $\pmb\mu_E = \tfrac{42}{3} =14,\ \pmb\mu_S=\tfrac{39}{3}=13$ and $\pmb\mu_\chi = (4+0-1)/3 = 1$.
\end{itemize}
Since vertex valencies and \piv\ status are unchanged by the grouping operation, the values of $\th$ and $\phi$ are unchanged from those found earlier in Section 5. So $\th = \tfrac73$ and $\phi= \tfrac16$ in both cases, A and B; thus $2\th/(\th-2) = 14$ and $2(\th-\phi)/(\th-2) = 13$.

Therefore, in case A, the formulae of (\ref{tag1}) yield $\pmb\mu_E = 14\,\pmb\mu_\chi = 14 \times \tfrac34 = \tfrac{42}{4}$ and $\pmb\mu_S = 13\, \pmb\mu_\chi = 13 \times \tfrac34 = \tfrac{39}{4}$, agreeing with the direct calculation above.

In case B, these formulae yield $\pmb\mu_E = 14\, \pmb\mu_\chi =  14$ and $\pmb\mu_S = 13\, \pmb\mu_\chi = 13$, agreeing with the direct calculation.

\section*{Acknowledgement}

Most of the ideas and analysis in this paper were presented in our technical report \cite{cots95} written in 1995. We feel that the results of that report, neglected by us for twenty years due to our differing career paths and changing interests, should now be placed online --- albeit with a clearer introduction and motivation based on geometric graphs. The results of the paper have, to our knowledge, not appeared elsewhere since 1995. In order to conform with our 1994 published paper on falling-leaf \tes s \cite{cots94} and to keep the historical setting of our 1995 technical report, with its connection to some other papers \cite{cow78, cow79, cow80}, notations and techniques of proof based on ergodic methods have not been altered greatly. A alternative theory based of Palm measures is, of course, possible --- and is being considered.

\bibliographystyle{plain}
\bibliography{CowanTsang_bibtex}   
\end{document}